\documentclass[12pt]{amsart}

\usepackage{amsmath,amssymb,amscd,amsfonts}

\newtheorem{thm}{Theorem}[section]

\newtheorem{prop}[thm]{Proposition}

\theoremstyle{remark}

\newcommand{\C}{\mathbb C}

\newcommand{\Q}{\mathbb Q}

\newcommand{\pp}{\mathbb P}
 \newcommand{\Alb}{\operatorname{Alb}}

\newcommand{\OO}{\mathcal O}

\newcommand{\fie}{\varphi}

\numberwithin{equation}{section}

\begin{document}
\title[On a theorem of Xiao Gang ]{
A note on a theorem of Xiao Gang }
\thanks{2000 Mathematics Subject Classification: 14J29}
\author{Margarida Mendes Lopes }
\date{}

\begin{abstract} In 1985  Xiao Gang proved that
the bicanonical system of    a  complex surface $S$
of general type with $p_2(S)>2$ is not composed of a pencil 
(\cite{xiaocan}).
In this note a new proof of this theorem is presented.
\end{abstract}

\maketitle

\section{Introduction}

In 1985 in \cite{xiaocan} Xiao Gang proved
that the bicanonical image of a (minimal) surface of
general type $S$ is not a surface (i.e., the bicanonical system is 
composed of a pencil) if and only if $K_S^2=1$,
$p_g(S)=0$, i.e., if and only if $p_2(S)=2$.  

When in the end of
the 80's it was finally proven that $| 2K_S|$ is base point
free, whenever $p_g\geq 1$, 
the part of this theorem concerning surfaces with $p_g\geq 1$ became 
trivial.

The aim of this note is to give a brief new  proof of this theorem of Xiao
Gang, using several results which by now are standard techniques of
surface theory.

\smallskip
\paragraph{\em Notations and conventions} A {\it surface} is an algebraic
projective surface  over $\C$. No distinction is made between line
bundles and divisors on a smooth variety, and  additive and multiplicative
notation are used interchangeably. Linear equivalence is denoted by 
$\equiv$ and
numerical equivalence by $\sim$. Given a linear system  $ | D| $  on a
surface, the corresponding rational map is denoted by $\fie_D$, and
$| D| $ is {\it composed of a pencil} if dim Im $\fie_D=1$. The remaining 
notation is standard in
algebraic geometry.
\smallskip
\paragraph{\em Acknowledgements}         This work was written     during 
a
stay
of the author
 at the institute ``Simon Stoilow'' of the Romanian Academy, under the
program Eurrommat.  Thanks     
are due
  to the Institute for the wonderful
hospitality and the good organization  found there.

\section{Preliminaries}
Here we list the results which will be needed in the sequel.

\begin{prop} [\cite{xiao}, Th\'eor\`eme 2.2]
\label{g2} Let $S$ be a minimal surface of general type with
$p_g(S)=0$. If $S$ has a genus $2$ fibration $f\colon S\to \pp^1$, then
$K_S^2\leq 2$.
\end{prop}

We will need also:
\begin{prop} [\cite{beau}, Corollaire 5.8]\label{etale}
A minimal surface of general type $S$ satisfying $q(S)=0$ and  $K_S^2\leq
2\chi(\OO_S)$ has no irregular   \'etale covers.    
\end{prop}

 \begin{thm}[De Franchis, \cite{defra},
 cf. \cite{cetraro}, \cite{topology}]
\label{deF}
 Let $S$ be a smooth surface with $p_g(S)=q(S)=0$ and $\pi\colon Y\to S$ a
smooth double cover with $q(Y)>0$. Then:
 \begin{enumerate}
 \item the Albanese map of\, $Y$ is a fibration $\alpha\colon Y\to B$ over 
a
curve $B\subset\Alb S$;
 \item there exist a fibration $g\colon S\to\pp^1$ and a degree $2$ map
$p\colon B\to\pp^1$ such that $p\circ \alpha =g\circ \pi$.
 \end{enumerate}
\end{thm}

Finally we recall: 
\begin{thm} \label{basepoints} The bicanonical system $|2K_S|$ of a
 minimal surface 
of general type $S$ 
is base point free, if $p_g(S)\geq 1$ or $K_S^2\geq 5$.
\end{thm}
This result is due to the work of several authors. For the corresponding 
references
see the  survey paper  \cite{ciro}.

\section {The proof of Xiao Gang's theorem}
 The theorem will follow from Theorem \ref{basepoints} and the following:

\begin{thm}\label{pencil}    Let $S$ be a minimal surface of general type 
with $p_g=0$. Then $|2K_S|$ is not composed of a
pencil
 if and only if $K_S^2> 1$.
\end{thm}

\begin{proof} Since $\chi(\OO_S)=1$, one has $h^0(S,2K_S)=K_S^2+1$ and 
therefore, if $K_S^2=1$, $|2K_S|$ is a pencil.

Suppose now that $K_S^2>1$ and that $|2K_S|$ is composed of a pencil with 
general fibre $F$. Since $q=0$, $|F|$ is a rational pencil and so
$2K_S\equiv dF+Z$, where $d=K_S^2$ and $Z$ is an effective divisor 
possibly zero.

Now, because $K_S$ is nef, we have $2K_S^2\equiv dK_SF+K_SZ\geq dK_SF$ and
therefore $K_SF\leq 2$. The index theorem  yields $K_S^2F^2\leq
(K_SF)^2$ where if equality holds  then for some $a,b\in \Q$, $aK_S\sim 
bF$. 
Also by the adjunction formula $K_SF$ and $F^2$ have
the same parity. By the assumption $K_S^2>1$, one obtains the
following  the numerical possibilities:

i) $K_SF=2$, $F^2=2$ and $K_S\sim F$;

ii) $K_SF=2$, $F^2=0$.

Now  case i) does not occur. Suppose
otherwise. Then  $2K_S\equiv 2F$ and so $\eta:=K_S-F$ is a
2-torsion divisor. The \'etale double cover $p:Y\to S$ associated to 
$\eta$ satisfies $\chi(\OO_Y)=2$ and
$p_g(Y)=h^0(S,K_S)+h^0(S,K_S+\eta)=h^0(S,K_S)+h^0(S,F)=2$. So $Y$ is
irregular. This is a contradiction to proposition 
\ref{etale},
 and case i) is excluded.

 For case ii) notice that, anyway, $K_S^2=2$, 
because in this case   the pencil $|F|$ is a genus $2$ fibration and 
therefore $K_S^2\leq 2$, by proposition \ref{g2}.

Then one has $2K_S\equiv 2F+Z$, where $Z>0$ is such that
$ K_SZ=0$, and so every irreducible curve in   $Z$ is a $-2-$curve.
The effective divisor $Z$ can be decomposed as   
$Z=2Z_0+Z_1$ where $Z_0, Z_1$ are  effective divisors, and $Z_1$ is
reduced.

 If $Z_1=0$, then   
 $\eta:=K_S-F-Z_0$ is a 2-torsion divisor and the same argument as above 
leads again to a contradiction. 

Suppose now that $Z_1\neq 0$.  Since $Z_1=2(K_S-F-Z_0)$,  $\theta Z_1$
is even for any irreducible component $\theta$ of $Z_1$.
On the other hand the dual graph of the configuration of curves in
$Z_1$ is a union of trees  and thus, because
$Z_1$ is reduced, 
 necessarily $Z_1$ is a disjoint union of $p$
irreducible  $-2-$curves. So $Z_1$  is a smooth divisor and we can 
consider                   the double cover $p:Y'\to S$ branched on $Z_1$ 
and defined
by the relation 
$Z_1=2(K_S-F-Z_0)$.

The standard double cover formulas yield

\vskip0.2truecm\noindent
$\chi ( \OO_{Y'})=2\chi(\OO_S)+\frac{1}{2}
(K_S-F-Z_0) (2K_S-F-Z_0)=2- \frac{p}{4}         $ ; 
\vskip0.2truecm\noindent$K_{Y'}^2=2(2K_S-F-Z_0)^2=4-p  $ ; 
\vskip0.2truecm\noindent$p_g(Y')=h^0(S, {\OO_S}(2K_S-F-Z_0))+h^0(S, 
{\OO_S}(K_S))= 2.$

\smallskip
Since the surface $Y'$ is of course of general type,
the only possibility is that $p=4$ ($\chi(\OO_{Y'})=1$). Now $p_g(Y')=2$ 
yields $q(Y')=2$ and
thus, by De Franchis theorem \ref{deF}, $Y'$ is not of Albanese general 
type.

On the other hand, the minimal model $Y$ of $Y'$ is obtained by
contracting the four exceptional curves which are the inverse images of 
the four components of $Z_1$ and as such satisfies $K_Y^2=4$. If we denote 
by
 $f$ the genus of a general fibre of the Albanese pencil of
$Y$, we  have then a contradiction to Arakelov's inequality
$K_Y^2\geq 8(f-1) (q(Y)-1)$ (see, e.g., \cite{bea}).

 So also this case does not occur and the theorem is proven.
\end{proof}

\smallskip\noindent {\bf Remark.} Very recently Meng Chen and E.Viehweg 
also found
a different proof of Xiao's theorem, which uses vanishing theo\-rems
for\- $\Q-$divisors (see \cite{MV}).

\bigskip 
\begin{tabbing} 1749-016 Lisboa, PORTUGALxxxxxxxxx \kill
Margarida Mendes Lopes\\ CMAF 
\\ Universidade de Lisboa  \\ Av. Prof. Gama Pinto, 2
\\ 1649-003 Lisboa, PORTUGAL \\ 
mmlopes@ptmat.ptmat.fc.ul.pt 
\end{tabbing}

\end{document}